\documentstyle[12pt]{article}

\newtheorem{theorem}{Theorem}
\newtheorem{lemma}{Lemma}
\newtheorem{definition}{Definition}

\def\diam{\mathop{\rm diam}}

\def\EXP{{\bf E}}
\def\PROB{{\bf P}}

\def\P{{\cal P}}

\def\Rd{{\cal R}^d}

\def\proof{\medskip  \par \noindent {\bf Proof.} \ \ }

\def\qedskip{\smallskip\noindent}
\def\qed{\hfill $\Box$ \qedskip}              %(end of proof after text)
\newcommand{\be}{\begin{equation}}
\newcommand{\ee}{\end{equation}}

\begin{document}
\bibliographystyle{plain}

\title{
A simple randomized algorithm for 
sequential prediction of ergodic time series\\
{\small (Appeared in:  IEEE Trans. Inform. Theory  45  (1999),  no. 7, 2642--2650.)}
}

\author{
L\'aszl\'o Gy\"orfi \\
\and
G\'abor Lugosi
\and
Guszt\'av Morvai \\
}

\date{}
\setcounter{page}{0}

\maketitle

\begin{abstract}
We present a simple randomized procedure for 
the prediction of a binary sequence. 
The algorithm uses ideas from recent
developments of the theory of the prediction
of individual sequences.
We show
that if the sequence is a realization of a stationary
and ergodic random process then the average
number of mistakes converges, almost surely, 
to that of the optimum, given by the Bayes 
predictor.
The desirable finite-sample properties of the predictor
are illustrated by its performance
 for Markov processes. In such cases the 
predictor exhibits near optimal behavior
even without knowing the order of the Markov
process. Prediction with side information
is also considered.
\end{abstract}
\newpage

\section{Introduction}

We address the problem of sequential prediction
of a binary sequence. 
A sequence of bits $y_1,y_2,\ldots \in \{0,1\}$ is
hidden from the predictor.
At each time instant $i=1,2,\ldots$,
the predictor is asked to guess
 the value of the next outcome
$y_i$ with knowledge of the past $y_1^{i-1}=(y_1,\ldots,y_{i-1})$
(where $y_1^0$ denotes the empty string). 
Thus, the predictor's decision, at time $i$, is based
on the value of $y_1^{i-1}$.
We also assume that the predictor has access to
a sequence of 
%LG
independent, identically distributed (i.i.d.) 
random variables $U_1,U_2,\ldots$,
uniformly distributed on $[0,1]$, so that the predictor
can use $U_i$ in forming a randomized decision for
$y_i$. Formally, the strategy of the predictor
is a sequence $g=\{g_i\}_{i=1}^{\infty}$ of decision functions
\[
   g_i: \{0,1\}^{i-1} \times [0,1] \to \{0,1\}
\]
and the randomized prediction formed at time $i$ is
$g_i(y_1^{i-1},U_i)$. The predictor pays a unit penalty
each time a mistake is made. After $n$ rounds of play,
the {\em normalized cumulative loss} 
on the string $y_1^n$ is
\[
  L_1^n(g,U_1^n) 
   = \frac{1}{n} \sum_{i=1}^n I_{\{g_i(y_1^{i-1},U_i) \neq y_i\}},
\]
where $I$ denotes the indicator function.
 When no confusion
is caused, or when the predictor does not randomize,
we will simply write $L_1^n(g)=L_1^n(g,U_1^n)$.
In general, we denote the average number of mistakes
between times $m$ and $n$ by
\[
  L_m^n(g,U_m^n) 
   = \frac{1}{n-m+1} \sum_{i=m}^n I_{\{g_i(y_1^{i-1},U_i) \neq y_i\}}.
\]
We also write
\[
\widehat{L}_1^n(g) = \EXP L_1^n(g,U_1^n) \quad \mbox{and} \quad
\widehat{L}_m^n(g) = \EXP L_m^n(g,U_m^n)
\]
for the expected loss of the randomized strategy $g$.
(Here the expectation is taken with respect to 
the randomization $U_1^n$.)

In this paper we assume that $y_1,y_2,\ldots$ are
realizations of the random variables $Y_1,Y_2,\ldots$
drawn from the binary-valued stationary and ergodic process
$\{Y_n\}_{-\infty}^{\infty}$. We assume that
the randomizing variables $U_1,U_2,\ldots$ are independent 
of the process $\{Y_n\}_{-\infty}^{\infty}$.

In this case there is a fundamental limit for
the predictability of the sequence. This is
stated in the next theorem whose proof may 
be found in Algoet \cite{Alg94}. 

\begin{theorem} (Algoet \cite{Alg94}) 
\label{bayes}
For any prediction strategy $g$ and stationary ergodic
process $\{Y_n\}_{-\infty}^{\infty}$,
\[
   \liminf_{n\to \infty} L_1^n(g) \ge L^* \quad \mbox{almost surely,}
\]
where 
\[
  L^* = \EXP\left[
\min\left(\PROB\{Y_0=1|Y_{-\infty}^{-1}\},\PROB\{Y_0=0|Y_{-\infty}^{-1}\}
  \right) \right]
\]
is the minimal (Bayes) probability of error of any decision
for the value of $Y_0$ based on the infinite past $Y_{-\infty}^{-1}=
(\ldots,Y_{-3},Y_{-2},Y_{-1})$.
\end{theorem}

\noindent
Based on Theorem \ref{bayes}, the following definition
is meaningful:

\begin{definition}
A prediction strategy $g$ is called {\em universal}
if for all stationary and ergodic processes $\{Y_n\}_{-\infty}^{\infty}$,
\[
 \lim_{n \to \infty} L_1^n(g) = L^* \quad \mbox{almost surely}.
\]
\end{definition}

Therefore, universal strategies asymptotically achieve
the best possible loss for all ergodic processes. 
The first question is, of course,
if such a strategy exists. The affirmative answer 
follows from a more general result of Algoet \cite{Alg94}.
Here we give an alternative proof which is based on
 earlier results of Ornstein and Bailey.

\begin{theorem} (Algoet \cite{Alg94})
\label{letezik}
There exists a universal prediction scheme.
\end{theorem}

\proof
Ornstein \cite{Orn78} proved that there exists a sequence
of functions $f_i: \{0,1\}^{i} \to [0,1]$, $i=1,2,\ldots$ 
  such that for all ergodic processes
$\{Y_n\}_{-\infty}^{\infty}$,
\be
\label{supposeconsistentestimator}
  \lim_{n\to \infty} f_n(Y_{-n}^{-1}) = \PROB\{Y_0=1|Y_{-\infty}^{-1}\}
   \quad \mbox{almost surely}.
\ee
Bailey \cite{Bai76} observed that for such estimators,
for all ergodic processes
\be
%LG indexek kijavitva ==> nezzetek meg, hogy jol csinaltam-e!
\label{bailey}
 \lim_{n\to \infty}  \frac{1}{n}\sum_{i=1}^n |(f_{i-1}(Y_1^{i-1})
 - \PROB\{Y_{i}=1|Y_{-\infty}^{i-1}\}| =0
   \quad \mbox{almost surely}.
\ee
Indeed, (\ref{supposeconsistentestimator}) and Breiman's generalized ergodic
theorem (see Lemma \ref{Breimanergod} in the Appendix) yield (\ref{bailey}).

Once such a sequence $\{f_i\}$ of estimators is available, we may 
define a (non-randomized) prediction scheme by the plug-in predictor
\[
 g_n(y_1^{n-1}) = \left\{ \begin{array}{ll}
1 & \mbox{if} \ \ f_{n-1}(y_1^{n-1}) \ge
\frac{1}{2} \\
0 & \mbox{otherwise.}
                  \end{array} \right.
\]
It is well-known that the probability of error of such a plug-in
predictor may be bounded by the $L_1$ error of the estimator
it is based on. In particular, by 
%LG 
a simple inequality appearing in the proof of \cite[Theorem~2.2]{DeGyLu95},
\[
\PROB\left\{g_n(Y^{n-1}_1)\neq Y_n|Y^{n-1}_{-\infty}\right\}-
\PROB\left\{g^*(Y^{n-1}_{-\infty})\neq Y_n|Y^{n-1}_{-\infty} \right\}\le
2\left|f_{n-1}(Y^{n-1}_1)-\PROB\left\{Y_n=1|Y^{n-1}_{-\infty}\right\}\right|.
\]
Therefore,
\begin{eqnarray*}
|L_1^n(g)-L^*|&\le& 
\left|L_1^n(g)-{1\over n} \sum_{i=1}^{n} 
\PROB\left\{g_i(Y^{i-1}_1)\neq Y_i|Y^{i-1}_1\right\}\right| \\*
& & + 
{1\over n} \sum_{i=1}^{n} 
\left|\PROB\left\{g_i(Y^{i-1}_1)\neq Y_i|Y^{i-1}_1\right\}-
\PROB\left\{g^*(Y^{i-1}_{-\infty})\neq
Y_{i}|Y^{i-1}_{-\infty}\right\}
 \right| \\*
& & +
\left|{1\over n} \sum_{i=1}^{n} 
\PROB\left\{g^*(Y^{i-1}_{-\infty})\neq
Y_{i}|Y^{i-1}_{-\infty}\right\}
 -L^*\right| \\
&\le& 
\left|L_1^n(g)-{1\over n} \sum_{i=1}^{n} 
 \PROB\left\{g_i(Y^{i-1}_1)\neq Y_i|Y^{i-1}_1\right\}\right| \\*
& & + 
{2\over n} \sum_{i=1}^{n} 
\left|f_{i-1}(Y^{i-1}_1)-\PROB\left\{Y_i=1|Y^{i-1}_{-\infty}\right\}\right|
\\*
& & +
\left|{1\over n} \sum_{i=1}^{n} 
\PROB\left\{g^*(Y^{i-1}_{-\infty})
 \neq Y_{i}|Y^{i-1}_{-\infty}\right\}-L^*\right|.
\end{eqnarray*}
The first term of the right-hand side tends to zero almost surely by
%LG the martingale convergence theorem.
the Hoeffding-Azuma inequality (Lemma \ref{azuma} in the Appendix)
and the Borel-Cantelli lemma.
The second one converges to zero almost surely by
(\ref{bailey}) and the third term tends to zero almost surely by the ergodic
theorem. 
\qed

It was Ornstein \cite{Orn78} who first proved the
existence of estimators satisfying (\ref{supposeconsistentestimator}).
This was later generalized by Algoet \cite{Alg92}.
A simpler
estimator with the same convergence property
was introduced by Morvai, Yakowitz, and Gy\"orfi
\cite{MoYaGy96}.
Unfortunately, even the simpler estimator
needs so
large amounts of data that its practical
use is unrealistic. By this we mean that 
even for ``simple'' i.i.d.\ or Markov processes
the rate of convergence of the estimator 
is very slow.
Motivated by the need of a practical estimator,
Morvai, Yakowitz, and Algoet \cite{MoYaAl97}
introduced an even simpler algorithm. However,
it is not known whether their estimator satisfies 
(\ref{supposeconsistentestimator}), and we do not even know whether
the corresponding predictor is universal.
The purpose of this paper is to introduce a new
simple universal predictor whose finite-sample
performance for Markov processes promise practical
applicability.

\newpage

\section{A simple universal algorithm}
\label{alg}

In this section we present a simple prediction
strategy, and prove its universality. It is motivated by
some recent developments from the theory of
the prediction of individual sequences (see, e.g.,
Vovk \cite{Vov90},
Feder, Merhav, and Gutman~\cite{FeMeGu92},
Littlestone and Warmuth \cite{LiWa94},
Cesa-Bianchi et al. \cite{CeFrHaHeScWa97}).
These methods predict according to a combination
of several predictors, the so-called {\em experts}.

The main idea in this paper is that if 
the sequence to predict is drawn from
a stationary and ergodic process, combining
the predictions of a small and simple set of 
appropriately chosen predictors (the so-called experts) suffices to
achieve universality.

First we define an infinite sequence
of experts $h^{(1)},h^{(2)},\ldots$
as follows: 
Fix a positive integer $k$, and for each $n\ge 1$, $s \in \{0,1\}^k$
and $y \in \{0,1\}$ 
define the function 
$\widehat{P}_n^k:\{0,1\} \times \{0,1\}^{n-1} \times \{0,1\}^k \to [0,1]$
by
\be
\label{freq}
{\widehat P}_n^k(y,y_1^{n-1},s)={
\left|\{k< i<n: y^{i-1}_{i-k}=s, y_i=y\}  \right|
\over 
\left|\{k< i < n: y^{i-1}_{i-k}=s\}\right| } \quad \mbox{for $n>k+1$} ,
\ee
where $0/0$ is defined to be $1/2$.
Also, for $n \le k+1$ we define ${\widehat P}_n^k(y,y_1^{n-1},s)=1/2$.
In other words, ${\widehat P}_n^k(y,y_1^{n-1},s)$ is the proportion
of the appearances of the bit $y$ following
the string $s$ among all appearances of $s$ in the sequence $y_1^{n-1}$.

\noindent
The expert $h^{(k)}$ is
a sequence of functions 
$h^{(k)}_n:\{0,1\}^{n-1}\to \{0,1\}$, $n=1,2,\ldots$
defined by
\[
  h^{(k)}_n(y_1^{n-1}) = \left\{ \begin{array}{ll}
            0 & \mbox{if} \ \  {\widehat 
P}_n^k(0,y_1^{n-1},y_{n-k}^{n-1})>{1\over 2} \\
            1 & \mbox{otherwise,} 
             \end{array}  \right. \quad n=1,2,\ldots .
\]
That is, expert $h^{(k)}$ is a (nonrandomized) prediction
strategy, which looks for all appearances
of the last seen string $y_{n-k}^{n-1}$ of length $k$
in the past and predicts according to the larger of the
relative frequencies of 0's and 1's following
the string. We may call $h^{(k)}$ a {\em $k$-th order
empirical Markov strategy.} 

The proposed prediction algorithm proceeds as follows:
Let $m=0,1,2,\ldots$ be a nonnegative integer.
For $2^m \le n < 2^{m+1}$, the prediction is based
upon a weighted majority of predictions of the experts
$h^{(1)},\ldots,h^{(2^{m+1})}$ as follows:
\[
  g_n(y_1^{n-1},u) = \left\{ \begin{array}{ll}
   0 & \mbox{if} \ \ 
\displaystyle{u> 
\frac{\sum_{k=1}^{2^{m+1}} 
h_n^{(k)}(y_1^{n-1}) w_n(k)}
     {\sum_{k=1}^{2^{m+1}} w_n(k)} }  \\
   1 & \mbox{otherwise,}
                     \end{array} \right. \quad n=1,2,\ldots,
\] 
where $w_n(k)$ is the weight of expert $h^{(k)}$ defined
by the past performance of $h^{(k)}$ as
\[
w_{2^m}(k)=1 \quad \mbox{and} \quad
%LG az exponens kijavitva (a kumulativ hiba kell oda, nem az atlagos)
 w_n(k) = e^{-\eta_m (n-2^m) L_{2^m}^{n-1}(h^{(k)})} \ \mbox{for} 
 \ 2^m < n < 2^{m+1},
\]
where $\eta_m=\sqrt{8\ln(2^{m+1})/2^m}$.
Recall that
\[
 L_{2^m}^{n-1}(h^{(k)}) 
= \frac{1}{n-2^m} \sum_{i=2^m}^{n-1} 
I_{\left\{ h^{(k)}_i(y_1^{i-1}) \neq y_i \right\} }
\]
is the average number of mistakes made by expert $h^{(k)}$
between times $2^m$ and $n-1$. The weight of each expert is therefore
exponentially decreasing with the number of its mistakes on
this part of the data.

\medskip
\noindent
{\bf Remarks.}
1. The above-mentioned estimator of 
Morvai, Yakowitz, and Algoet \cite{MoYaAl97}
selects a value of $k$ in a certain data-dependent manner,
and uses the corresponding estimate $\widehat{P}_n^k$.
The new estimate, however, takes a mixture (weighted average)
of all possible values of $k$, with exponential 
weights depending on the past performance of
each component estimator. As Lemma \ref{expert}
below suggests, this technique guarantees a number
of errors almost as small as that of the {\em best}
expert (i.e., best value of $k$).

\medskip 
\noindent
2. Ryabko \cite{Rya88} proposed an estimator somewhat similar
in spirit to the predictor defined here. Ryabko
used a mixture of empirical Markov predictors,
and proved its universality for all stationary
and ergodic processes in a sense related
to the Kullback-Leibler divergence.
The idea of diversifying Markov strategies
also appears in Algoet \cite{Alg92}.

\medskip 
\noindent
3. Each time $n$ equals a power of two, all weights
are reset to 1, and a simple majority vote
is taken among the experts. This is necessary to
make the algorithm sequential and to be able to
incorporate more and more experts in the decision.
If the total length of the sequence to be predicted
was finite (say $n$) and known in advance, then
no such resetting would be necessary,
one could just use the first $n$ experts
as Lemma \ref{expert} below describes.
However, to achieve universality, an infinite
class of experts is necessary. As the first part of
the proof of Theorem \ref{cons} below shows,
we do not loose much by such a resetting of the weights.

\medskip 
\noindent
4. Related prediction schemes have been proposed
by Feder, Merhav, and Gutman~\cite{FeMeGu92} for
individual sequences. Their computationally
quite simple methods are shown to predict
asymptotically as well
as any finite-state predictor.

\medskip
\noindent
The main result of this section
is the universality of this simple prediction
scheme:

\begin{theorem}
\label{cons}
The prediction scheme $g$ defined above is universal.
\end{theorem}

\noindent
In the proof we use a beautiful result of
Cesa-Bianchi et al. \cite{CeFrHaHeScWa97}.
It states that, given a set of 
%LG kicsereltem a szakertok szamat N-rol K-ra, ez tenyleg jobb jeloles 
$K$ experts,
and a sequence of fixed length $n$, there
exists a randomized predictor whose number of mistakes
is not more than that of the best predictor
plus  $\sqrt{(n/2)\ln K}$ for {\em all
possible} sequences $y_1^n$. The simpler algorithm
and statement cited below is due to Cesa-Bianchi \cite{Ces97}:

\begin{lemma}
\label{expert}
Let $\tilde{h}^{(1)},\ldots,\tilde{h}^{(K)}$ be a finite collection
of prediction strategies (experts), and let $\eta >0$. 
Then if the prediction strategy $\tilde{g}$ is defined
by
\[
%LG a t-ket kicsereltem i-re
  \tilde{g}_i(y_1^{i-1},u) = \left\{ \begin{array}{ll}
   0 & \mbox{if} \ \ 
\displaystyle{u> 
\frac{\sum_{k=1}^{K}
\PROB{\left\{\tilde{h}^{(k)}(y_1^{i-1},U_i)=1\right\}}
 \tilde{w}_i(k)}
     {\sum_{k=1}^{K} \tilde{w}_i(k)} }  \\
   1 & \mbox{otherwise,}
                     \end{array} \right.
\]
$i=1,2,\ldots$, where  for all $k=1,\dots,K$
\[
\tilde{w}_1(k)=1 \quad \mbox{and} \quad
\tilde{w}_i(k) = e^{-\eta (i-1)\widehat{L}_1^{i-1}(\tilde{h}^{(k)})}, \ \
i>1
\]
then for every $n \ge 1$ and $y_1^n \in \{0,1\}^n$,
\[
\widehat{L}_1^n(\tilde{g}) \le \min_{k=1,\ldots,K}
 \widehat{L}_1^n(\tilde{h}^{(k)}) + \frac{\ln K}{\eta n} + \frac{\eta}{8}.
\]
In particular, if $N$ is a positive integer, and 
$\eta=\sqrt{8N^{-1}\ln K}$, then
\[
\widehat{L}_1^n(\tilde{g}) \le \min_{k=1,\ldots,K}
 \widehat{L}_1^n(\tilde{h}^{(k)}) + \frac{\sqrt{N}}{n}\sqrt{\frac{\ln K}{2}},
%LG
  \quad n \le N.
\]
\end{lemma}

\medskip
\noindent
{\bf Proof of Theorem \ref{cons}.}
%LG
Taking $K=2^{m+1}$ and $N=2^m$ in 
 Lemma \ref{expert}, we have that the expected number of errors 
committed by $g$ on a segment $2^m,\ldots,2^{m+1}-1$
is bounded, for any $y_{2^m}^{2^{m+1}-1}\in \{0,1\}^{2^m}$, as
\begin{eqnarray*}
\widehat{L}_{2^m}^{2^{m+1}-1}(g)
& = & \EXP\left[\frac{1}{2^m} \sum_{i=2^m}^{2^{m+1}-1} 
 I_{\{g_i(y_1^{i-1},U_i) \neq y_i\}} \right] \\*
 & \le & 
\min_{k \le 2^{m+1}} L_{2^m}^{2^{m+1}-1}(h^{(k)}) +
\sqrt{\frac{\ln (2^{m+1})}{2\cdot 2^m}} \\*
& = &
\min_{k=1,2,\ldots} L_{2^m}^{2^{m+1}-1}(h^{(k)}) +
\sqrt{\frac{\ln (2^{m+1})}{2\cdot 2^m}},
\end{eqnarray*}
where the last equality follows from the fact that 
for all 
$i < 2^{m+1}$, all
experts $h^{(k)}$ with $k\ge 2^{m+1}$
predict identically to $h^{(2^{m+1})}$.
(Note that since the predictors $h^{(k)}$ are deterministic,
for every $m$, 
$\widehat{L}_{2^m}^{2^{m+1}-1}(h^{(k)}) = L_{2^m}^{2^{m+1}-1}(h^{(k)})$.)

%LG tenyleg meg lehet javitani, de csak annak az
%LG aran, ha az n helyett az \overline{n} szerepel
%LG a korlatban, ugyhogy ezen nem valtoztattam.
%LG Masreszt viszont a \log_2 n+1 helyett \log_2 n
%LG eleg. Itt van a javitott levezetes.
Similarly, denoting $\overline{n}=2^{\lfloor \log_2 n\rfloor +1}$,
and invoking Lemma \ref{expert} with $K=\overline{n}$ and
$N=\overline{n}/2$,
\[
L_{{\overline{n}}/2}^{n}(g) \le 
\left(\min_{k=1,2,\ldots} L_{{\overline{n}}/2}^{n}(h^{(k)}) +
\frac{\sqrt{\overline{n}/2}}{n-{\overline{n}}/2+1}
\sqrt{\frac{\ln ({\overline{n}})}{2}} \right). 
\]
Therefore, 
for any sequence $y_1,y_2,\ldots$,
\begin{eqnarray*}
n\widehat{L}_1^n(g) 
& = &  \sum_{m=0}^{\lfloor \log_2 n \rfloor-1}
 2^m L_{2^m}^{2^{m+1}-1}(g)
 + (n-{\overline{n}}/2+1) L_{{\overline{n}}/2}^{n}(g) \\*
& \le &
 \sum_{m=0}^{\lfloor \log_2 n \rfloor-1}
 2^m \left(\min_{k=1,2,\ldots} L_{2^m}^{2^{m+1}-1}(h^{(k)}) +
\sqrt{\frac{\ln (2^{m+1})}{2\cdot 2^m}} \right)  \\*
& & + (n-{\overline{n}}/2+1)
\left(\min_{k=1,2,\ldots} L_{{\overline{n}}/2}^{n}(h^{(k)}) +
\frac{\sqrt{\overline{n}/2}}{n-{\overline{n}}/2+1}
\sqrt{\frac{\ln ({\overline{n}})}{2}} \right)  \\
& \le &
n \min_{k=1,2,\ldots} L_1^n(h^{(k)})
+
 \sum_{m=0}^{\lfloor \log_2 n \rfloor-1}
 \sqrt{\frac{2^m\ln (2^{m+1})}{2}} 
+ \sqrt{\frac{\overline{n}/2\ln \overline{n}}{2}}  \\*
&  = & 
n \min_{k=1,2,\ldots} L_1^n(h^{(k)})
+
 \sum_{m=0}^{\lfloor \log_2 n \rfloor}
 \sqrt{\frac{2^{m}\ln (2^{m+1})}{2}}.
\end{eqnarray*}
%G itt kezdodik
Denoting
$\mu=\lfloor \log_2 n\rfloor$, we may write
\begin{eqnarray*}
\sum_{m=0}^\mu\sqrt{2^m\ln(2^{m+1})\over 2}
&\le& \sqrt{\ln (2^{\mu+1})\over 2}\sum_{m=0}^\mu 2^{m/2}\\
&<& \sqrt{\ln(2^{\mu+1})\over 2}\,\cdot\,{2^{(\mu+1)/2}\over 
\sqrt{2}-1}\\
&=& c\sqrt{\overline{n}\log_2 \overline{n}\over 2},
\end{eqnarray*}
where
$$c={\sqrt{\ln 2}\over \sqrt{2}-1} \approx 2.01.$$
Thus, we obtain
\begin{eqnarray*}
\widehat{L}_1^n(g)
& \le & \min_{k=1,2,\ldots} L_1^n(h^{(k)})+{c\over n}\sqrt{\overline{n}
\log_2\overline{n}\over 2}\\
& \le & \min_{k=1,2,\ldots} L_1^n(h^{(k)})+c\sqrt{\log_2 n+1\over n}.
\end{eqnarray*}
%G eddig.

Noting that for any fixed sequence $y_1^n$,  $L_1^n(g,U_1^n)$
is a sum of $[0,1]$-valued independent random variables
whose expectation is $\widehat{L}_1^n(g)$, we may
use Hoeffding's inequality \cite{Hoe63} to see that
for any sequence $y_1^n$, and $\epsilon>0$,
\begin{equation}
\label{conc}
  \PROB\left\{ \left|L_1^n(g,U_1^n)-\widehat{L}_1^n(g)
         \right| > \epsilon \right\} \le 2e^{-2n\epsilon^2}.
\end{equation}
Therefore, if $L$ is  now evaluated on the
random sequence $Y_1,Y_2,\ldots$, we obtain
\begin{eqnarray*}
\limsup_{n\to \infty} L_1^n(g,U_1^n)
& \le & \limsup_{n\to \infty}
\left( \min_{k=1,2,\ldots}
 L_1^n(h^{(k)}) + 
c \sqrt{\frac{\log_2 n + 1}{n}} \right).
                                \\*
& = &
\limsup_{n\to \infty} 
\min_{k=1,2,\ldots}
 L_1^n(h^{(k)})
 \quad \mbox{almost surely}.
\end{eqnarray*}
Thus, it remains to show that for any ergodic process $Y_1,Y_2,\ldots$,
\be
\label{eq0}
\limsup_{n\to \infty}
\min_{k=1,2,\ldots}
 L_1^n(h^{(k)}) \le L^*  \quad \mbox{almost surely}.
\ee
This will follow easily from the
following lemma:

\begin{lemma}
\label{fststat}
For any $k \ge 1$,
\[
\limsup_{n\to\infty} L_1^n(h^{(k)})\le 
L^* + \epsilon_k
 \quad \mbox{almost surely},
\]
where $\epsilon_k>0$ is such that $\lim_{k\to \infty} \epsilon_k=0$.
\end{lemma}

\medskip
\noindent
{\bf Remark.} If the process $\{Y_n\}$ happens to be $m$-th order
Markov, then it is easy to see that $\epsilon_k=0$ for
all $k\ge m$. The performance of the predictor for
such processes is investigated in the next section.

\proof
Introduce 
$$
{\tilde L}_1^n (h^{(k)})
= {1\over n} \sum_{i=1}^{n} 
\PROB\{Y_i\neq h_i^{(k)}(Y^{i-1}_{1})|Y_{-\infty}^{i-1}\}.
$$
By Lemma \ref{azuma} in the Appendix
we immediately obtain
$$
\lim_{n\to\infty} \left| L_1^n(h^{(k)})-{\tilde L}_1^n(h^{(k)}) \right| = 0
\quad \mbox{almost surely. }
$$
Therefore, it suffices to show that
$\limsup_{n\to\infty} \tilde{L}_1^n(h^{(k)})\le 
L^* + \epsilon_k$ almost surely. To this end,
first we study the asymptotic behavior of the quantity
$\PROB\{Y_0\neq h_n^{(k)}(Y^{-1}_{-n+1})|Y_{-\infty}^{-1}\}$.
Notice that 
\begin{eqnarray}
\label{terms}
\PROB\{Y_0\neq h_n^{(k)}(Y^{-1}_{-n+1})|Y_{-\infty}^{-1}\}
&\le&
I_{A}\PROB\{Y_0\neq h_n^{(k)}(Y^{-1}_{-n+1})|Y_{-\infty}^{-1}\}
\nonumber \\
& & + I_{B_k}\PROB\{Y_0\neq h_n^{(k)}(Y^{-1}_{-n+1})|Y_{-\infty}^{-1}\} 
\nonumber \\
& & + I_{C_k}\PROB\{Y_0\neq h_n^{(k)}(Y^{-1}_{-n+1})|Y_{-\infty}^{-1}\} 
\nonumber \\
& & + I_{D_k^c}\PROB\{Y_0\neq h_n^{(k)}(Y^{-1}_{-n+1})|Y_{-\infty}^{-1}\}
\end{eqnarray}
where 
\[
A=\left\{ \PROB\{Y_0=1|Y_{-\infty}^{-1}\}={1\over 2}\right\},
\]
\[
B_k=\left\{\PROB\{Y_0=1|Y_{-\infty}^{-1}\}<
{1\over 2} \ \mbox{and} \ \PROB\{Y_0=1|Y_{-k}^{-1}\}<{1\over 2}\right\},
\]
\[
C_k=\left\{\PROB\{Y_0=0|Y_{-\infty}^{-1}\}<{1\over 2}
 \ \mbox{and} \ \PROB(Y_0=0|Y_{-k}^{-1})<{1\over 2}\right\},
\]
and $D_k=A\cup B_k\cup C_k$.
Notice that 
\begin{eqnarray}
\label{conderrorterm}
\lefteqn{
\PROB\{Y_0\neq h_n^{(k)}(Y^{-1}_{-n+1})|Y_{-\infty}^{-1}\}  }
\nonumber \\*
& & = 
\PROB\{Y_0=1|Y_{-\infty}^{-1}\}
I_{\{{\widehat P}_n^k(1,Y^{-1}_{-n+1}|Y^{-1}_{-k})\le {1\over 2}\}}
+
\PROB\{Y_0=0|Y_{-\infty}^{-1}\}
I_{\{{\widehat P}_n^k(1,Y^{-1}_{-n+1}|Y^{-1}_{-k})> {1\over 2}\}}. 
\end{eqnarray}
Now we examine the four terms on the right-hand side of (\ref{terms}).
For the first term 
(\ref{conderrorterm}) yields
\begin{eqnarray*} 
I_{A}\PROB\{Y_0\neq 
h_n^{(k)}(Y^{-1}_{-n+1})|Y_{-\infty}^{-1}\}
&=& 
I_{A}{1\over 2} \\
&=&
I_{A}\min\left(\PROB\{Y_0=1|Y_{-\infty}^{-1}\},
               \PROB\{Y_0=0|Y_{-\infty}^{-1}\} \right).
\end{eqnarray*}
For the second term observe that under $B_k$, 
for sufficiently large $n$,
\[
{\widehat P}_n^k(1,Y_{-n+1}^{-1},Y_{-k}^{-1})<{1\over 2}
\quad \mbox{ almost surely,}
\]
and therefore by (\ref{conderrorterm}) we have
\[
\lim_{n\to\infty} I_{B_k}\PROB\{Y_0\neq
h_n^{(k)}(Y^{-1}_{-n+1})|Y_{-\infty}^{-1})
=I_{B_k}\min\left(\PROB\{Y_0=1|Y_{-\infty}^{-1}\},
                  \PROB\{Y_0=0|Y_{-\infty}^{-1}\} \right)
 \quad \mbox{a.s.}
\]
For the third term we obtain similarly
\[
\lim_{n\to\infty} I_{C_k}\PROB\{Y_0\neq
h_n^{(k)}(Y^{-1}_{-n+1})|Y_{-\infty}^{-1})
=I_{C_k}\min\left(\PROB\{Y_0=1|Y_{-\infty}^{-1}\},
                  \PROB\{Y_0=0|Y_{-\infty}^{-1}\} \right)
 \quad \mbox{a.s.}
\]
The last term is simply bounded by
\[
I_{D_k^c} \PROB\{Y_0\neq h_n^{(k)}(Y^{-1}_{-n+1})|Y_{-\infty}^{-1}\} \le
I_{D_k^c}.
\]

\noindent
Combining all these bounds, we obtain
\begin{equation}
\label{uppbound}
\PROB\{Y_0\neq h_n^{(k)}(Y^{-1}_{-n+1})|Y_{-\infty}^{-1}\} 
\le
I_{D_k}\PROB\{Y_0\neq h_n^{(k)}(Y^{-1}_{-n+1})|Y_{-\infty}^{-1}\}+I_{D_k^c}
\end{equation}
and 
\begin{equation}
\label{limit}
\lim_{n\to\infty}
I_{D_k}\PROB\{Y_0\neq h_n^{(k)}(Y^{-1}_{-n+1})|Y_{-\infty}^{-1}\}=
I_{D_k}\min\left(\PROB\{Y_0=1|Y_{-\infty}^{-1}\},
                 \PROB\{Y_0=0|Y_{-\infty}^{-1}\} \right)
 \quad \mbox{a.s.}
\end{equation}
 From (\ref{uppbound}) it is immediate that 
\[
\limsup_{n\to\infty} \tilde{L}_1^n(h^{(k)}) \le 
\lim_{n\to\infty}
{1\over n} \sum_{i=1}^{n} 
 I_{D_k}(T^i Y^{\infty}_{-\infty})
\PROB\{Y_i\neq h_i^{(k)}(Y^{i-1}_{1})|Y_{-\infty}^{i-1}\}+
\lim_{n\to\infty}
{1\over n} \sum_{i=1}^{n} 
I_{D_k^c}(T^i Y^{\infty}_{-\infty}),
\]
where $T$ denotes the left shift operator defined on doubly infinite
binary sequences $y_{-\infty}^{\infty} \in \{0,1\}_{-\infty}^{\infty}$.
By this inequality and (\ref{limit}), Breiman's generalized ergodic
theorem (see Lemma \ref{Breimanergod}
in the Appendix) implies 
\begin{eqnarray*}
\limsup_{n\to\infty} \tilde{L}_1^n(h^{(k)}) 
&\le&
\EXP\left[\min\left(\PROB\{Y_0=1|Y_{-\infty}^{-1}\},
              \PROB\{Y_0=0|Y_{-\infty}^{-1}\} \right) \right]
+ \PROB\{D_k^c\} \\*
&=& L^* + \PROB\{D_k^c\} \quad \mbox{almost surely.}
\end{eqnarray*}
Since by the martingale convergence theorem
$$
\lim_{k\to\infty}\PROB\{Y_0=1|Y^{-1}_{-k}\}
 =\PROB\{Y_0=1|Y^{-1}_{-\infty}\} \quad \mbox{almost surely,} 
$$
we have
$$\lim_{k\to\infty}\PROB(D_k^c)=0.$$
Taking
$\epsilon_k=\PROB\{D_k^c\}$, the proof of the lemma is complete. 
\qed

Now we return to the proof of Theoren~\ref{cons}.
\noindent
By  Lemma \ref{fststat}, 
for arbitrary $K$, 
\begin{eqnarray*}
\limsup_{n\to\infty} \min_{k=1,2,\ldots}
L_1^n(h^{(k)})
&\le&
\limsup_{n\to\infty} L_1^n(h^{(K)}) \\*
&\le & 
L^* + \epsilon_K.
\end{eqnarray*}
Since $K$ is arbitrary and $\epsilon_K\to 0$, (\ref{eq0}) is
established,
and the proof of the theorem is finished. \qed

\medskip
\noindent
{\bf Remarks.}
1. The proposed estimate is clearly easy to compute.
One merely has to keep track of 
the expected cumulative losses $L_{2^m}^{n-1}(h^{(k)})$
for $k=1,2,\ldots,n$. However, for large $n$, storing the
entire data history may be problematic. In such cases,
more efficient tree-based data structures, such as the
ones described by Feder, Merhav, and Gutman \cite{FeMeGu92},
may be applied. We do not investigate this issue further here.

\medskip
\noindent
2. We see from the analysis that for {\em any} sequence
$y_1,y_2,\ldots$ and for all $n$, 
\[
\widehat{L}_1^n\left(g\right) \le \min_{k=1,2,\ldots} 
L_1^n\left(h^{(k)}\right)
%G apro javitas 
+ 2.01\sqrt{\frac{\log_2 n+1}{n}},
\]
and that 
%LG
the difference $|L_1^n(g,U_1^n)-\widehat{L}_1^n\left(g\right)|$
between the actual loss and the expected loss  is 
$O_p(n^{-1/2})$.
%LG elorehoztam az O_p definiciojat, es korrigaltam
(For a sequence of random variables $\{X_n\}$
and sequence of nonnegative numbers $\{a_n\}$
we say that $X_n = O_p(a_n)$ if for every $\epsilon>0$
there exists a constant $c>0$ such that
$\limsup_{n \to \infty} \PROB\{|X_n| \ge ca_n\} < \epsilon$.)
%LG In other words, the algorithm is guaranteed to perform almost
%LG as well as the best experts.
The rate of convergence to $L^*$ depends on the behavior
of the best expert
%LG
for the time segment up to $n$. For example, in the next section we show
that for $m$-th order Markov processes the $m$-th expert
predicts very well, and this fact will suffice to derive
performance bounds for the proposed predictor.

\medskip 
\noindent
3. The proposed predictor is by no means the only
possibility. Different sets of experts
may be combined in a similar fashion,
and universality only depends on the behavior
of the best expert. If some additional information
is known (or suspected) about the process
to be predicted, this information may be
built in the definition of the experts.
We chose the empirical Markov strategies as
experts for convenience, and as we'll see it in the next section,
this choice pays off whenever the process
happens to be finite order Markov.

\newpage

\section{Markov processes}
\label{mark}

In this section we assume that the process 
to predict $\{Y_n\}_{-\infty}^{\infty}$ is 
(in addition to being stationary and ergodic)
$m$-th order Markov, that is, for any binary sequence
$y_{-\infty}^{-1}=(\ldots,y_{-2},y_{-1})$,
\[
  \PROB\{Y_0=1|Y_{-\infty}^{-1}=y_{-\infty}^{-1}\}
 =   \PROB\{Y_0=1|Y_{-m}^{-1}=y_{-m}^{-1}\},
\]
where $m$ is a positive integer.
 We show that the proposed
predictor achieves a nearly optimal performance
for any $m$ and for any such process,
even though the predictor does not use
the knowledge that the process is $m$-th order
Markov. The intuitive reason for such a behavior
is the following:
we have seen it in the previous
section that for {\em any} sequence,
\[
 L_1^n(g,U_1^n) \le \min_{k=1,2,\ldots} 
L_1^n\left(h^{(k)}\right)
+ 3\sqrt{\frac{\log_2 n+1}{n}}
+ O_p\left(\sqrt{\frac{1}{n}}\right).
\]
On the other hand, if the sequence is $m$-th
 order Markov, then there exists an expert,
namely $h^{(m)}$ with very good performance.

In order to simplify our analysis, we modify
the experts somewhat. They are defined as
before but the probability estimates of
(\ref{freq}) are now replaced by
\begin{equation}
\label{laplace}
\bar{P}_n^k(y,y_1^{n-1},s)={
\left|\{k< i<n: y^{i-1}_{i-k}=s, y_i=y\}  \right| +1
\over 
\left|\{k< i < n: y^{i-1}_{i-k}=s\}\right| +2}.
\end{equation}
In other words, the simple empirical frequency
counts are now replaced by the corresponding 
Laplace estimates. It is easy to see that 
all results of Section \ref{alg} remain
valid for the modified predictor.

\medskip
\noindent
{\bf Remark.}
The reason for this modification is that
this way we can appeal to a result of Rissanen \cite{Ris86}
which simplifies our analysis. We believe that
similar performance bounds are true for the original
predictor of Section \ref{alg}.

\medskip
In the next theorem we compare the performance
of our predictor to the universal lower bound $L^*$.
The statement only gives information about the expected
loss, but we believe this result already illustrates
the good behavior of the proposed predictor for Markov processes.

\begin{theorem}
\label{markov}
If the process to be predicted is a stationary and
ergodic $m$-th order Markov process, then the cumulative
loss $L_1^n(g)=L_1^n(g,U_1^n)$ of the
prediction
strategy of Section \ref{alg} (with the modified
estimates of (\ref{laplace}))
satisfies
\[
\EXP L_1^n(g) \le L^* + 2\sqrt{\frac{2^{m-1}\log n}{n}} + 
 3\sqrt{\frac{\log_2 n+1}{n}} 
   + \sqrt{\frac{c}{n}},
\]
where $c>0$ is a universal constant.
\end{theorem}

\proof
First note that (\ref{conc}) implies
\[
\EXP \left[| L_1^n(g,U_1^n) - \widehat{L}_1^n(g) ||Y_{-\infty}^{\infty}\right]
\le \int_0^{\infty} 2e^{-2n\epsilon^2} d\epsilon
\le \sqrt{\frac{\ln(2e)}{2n}}
\]
(see, e.g., \cite[page 208]{DeGyLu95}),
and therefore it suffices to investigate $\widehat{L}_1^n(g)$.
Recall also from the proof of Theorem \ref{cons} that 
for any input sequence,
\[
\widehat{L}_1^n(g) \le \min_{k=1,2,\ldots} L_1^n(h^{(k)})
 +  3\sqrt{\frac{\log_2 n+1}{n}}, 
\]
and, in particular,
\[
\widehat{L}_1^n(g) \le L_1^n(h^{(m)})
 +  3\sqrt{\frac{\log_2 n+1}{n}}.
\]
Thus, it suffices to show that for $m$-th order Markov
processes the performance of the $m$-th expert $h^{(m)}$
satisfies
\[
\EXP L_1^n(h^{(m)}) \le L^* + 2\sqrt{\frac{2^{m-1}\log n}{n}}
  + \sqrt{\frac{c}{n}}
\]
for some constant $c$.
To this end, observe that, on the one hand,
\begin{eqnarray*}
\EXP L_1^n(h^{(m)}) & = &
\EXP \left[ 
\frac{1}{n} \sum_{i=1}^n I_{\{h_i^{(m)}(Y_1^{i-1}) \neq Y_i\}} \right] \\*
& = & 
\EXP \left[ 
\frac{1}{n} \sum_{i=1}^n \PROB\{h_i^{(m)}(Y_1^{i-1}) \neq Y_i
    |Y_{-\infty}^{i-1}\} \right],
\end{eqnarray*} 
and on the other hand, by the Markov property,
\begin{eqnarray*}
  L^* & = & \EXP \left[ 
\frac{1}{n} \sum_{i=1}^n
 \min \left(\PROB\left\{Y_i=1|Y_{i-m}^{i-1}\right\},
             \PROB\left\{Y_i=1|Y_{i-m}^{i-1}\right\} \right) \right] \\*
& = & \EXP \left[
\frac{1}{n} \sum_{i=1}^n 
  \PROB\left\{ h^{(m,*)}(Y_{i-m}^{i-1}) \neq Y_i| Y_{-\infty}^{i-1} \right\} \right],
\end{eqnarray*}
where $h^{(m,*)}$ is the Bayes decision, given,
for any $s \in \{0,1\}^m$, by
\[
  h^{(m,*)}(s) = \left\{ \begin{array}{ll}
                  1 & \mbox{if $\PROB\{Y_0=1|Y_{-m}^{-1}=s\} \ge 1/2$} \\
                  0 & \mbox{otherwise.}
                 \end {array} \right.
\]  
(Note that the optimal predictor,
that is, the one which minimizes the probability
of error at every step predicts according to $h^{(m,*)}$.)

The above equalities imply that
\begin{eqnarray*}
\EXP L_1^n(h^{(m)}) - L^* & \le &
\frac{1}{n} \sum_{i=1}^n \EXP \left|
\PROB\left\{ h_i^{(m)}(Y_1^{i-1}) \neq Y_i |Y_{-\infty}^{i-1}\right\}  - 
 \PROB\left\{ h^{(m,*)}(Y_{i-m}^{i-1}) \neq Y_i| Y_{-\infty}^{i-1} \right\} 
 \right|
\\*
& \le &
\frac{2}{n} \sum_{i=1}^n \EXP \left|
\bar{P}_i^m(1,Y_1^{i-1},Y_{i-m}^{i-1})
-  \PROB\left\{ Y_i=1| Y_{1}^{i-1} \right\}  \right|, 
\end{eqnarray*}
where the second inequality follows by \cite[Theorem 2.2]{DeGyLu95}.
In the rest of the proof we simply apply some 
known results from the theory of universal prediction.
First, by applications of Jensen's and Pinsker's
inequalities (see Merhav and Feder \cite[eq. (20)]{MeFe98})
we obtain
\begin{eqnarray*}
\lefteqn{
\frac{2}{n} \sum_{i=1}^n \EXP \left|
\bar{P}_i^m(1,Y_1^{i-1},Y_{i-m}^{i-1})
-  \PROB\left\{ Y_i=1| Y_{1}^{i-1} \right\}  \right|   }  \\*
& \le &
2 \sqrt{\frac{1}{n} \sum_{i=1}^n \sum_{y_1^{i-1} \in \{0,1\}^{i-1}}
  \PROB\{Y_1^{i-1}=y_1^{i-1}\} \sum_{j=0}^1
 \PROB\{Y_i=j|Y_1^{i-1}=y_1^{i-1}\} 
 \log \frac{\PROB\{Y_i=j|Y_1^{i-1}=y_1^{i-1}\}}
           {\bar{P}_i^m(j,y_1^{i-1},y_{i-m}^{i-1})}}.
\end{eqnarray*}
Observe that on the right-hand side, under the square root
sign, we have the normalized Kullback-Leibler
divergence between the probability measure
of $Y_1^n$ and its estimate constructed as a
product of the Laplace estimates (\ref{laplace}).
But this divergence, for $m$-th order
Markov sources, is well-known to be bounded
by
\[
  \frac{2^m}{2n} \log n + O\left(\frac{1}{n}\right),
\]
see Rissanen \cite{Ris86}. This concludes the proof.
\qed

\medskip
\noindent
{\bf Remarks.}
1.
As Theorem \ref{markov} shows, by exponential weighting
of the empirical Markov strategies, 
the predictor automatically adapts to 
the unknown Markov order. Similar results,
though in different setup, are achieved
by Modha and Masry \cite{MoMa96},\cite{MoMa98} by 
complexity regularization.

\medskip
\noindent
2. Merhav, Feder, and Gutman, \cite{MeFeGu93} showed that
if the process is $m$-th order Markov, then
the randomized predictor $\tilde{h}^{(m)}$ defined by
\[
 \tilde{h}^{(m)}_i(y_1^{i-1},U) = \left\{ \begin{array}{ll}
 0 & \mbox{if} \ \  \widehat{P}_n^m(0,y_1^{n-1},y_{n-m}^{n-1})>{1\over 2} \\
 1 & \mbox{if} \ \  \widehat{P}_n^m(0,y_1^{n-1},y_{n-m}^{n-1})<{1\over 2} \\
 I_{\{U \ge 1/2\}} & \mbox{otherwise} 
                  \end{array} \right.
\]
achieves $\EXP L_1^n(\tilde{h}^{(m)}) - L^* \le C/n$,
where $C$ is a constant depending of the distribution
of the process. However, in an interesting contrast,
the best {\em distribution-free} upper bound for all
$m$-th order Markov processes is of the order of
$n^{-1/2}$. 
To illustrate this,
 consider the case $m=0$, that is, when
$\{Y_n\}$ is an i.i.d.\ process  with $\PROB\{Y_1=1\} = 1/2+\theta$,
and the predictor $\tilde{h}^{(0)}$ is based on
a majority vote of the bits appeared in the past.
In this case, for every $n$,
\[
 \sup_{\theta \in [-1/2,1/2]} 
 \left( \EXP L_1^n(\tilde{h}^{(0)}) - L^* \right) \ge c_1n^{-1/2},
\]
where $c_1$ is a universal constant.
%LG itt a bizonyitas vazlata
(This is straightforward to see by considering
$\theta = cn^{-1/2}$ for some small constant $c$,
and writing
\begin{eqnarray*}
\EXP L_1^n(\tilde{h}^{(0)}) - L^* 
& = &
\frac{1}{n} \sum_{i=1}^n \left[ 
\PROB\{ \tilde{h}^{(0)}(Y_1^{i-1},U_i) \neq Y_i\} 
 - \left(\frac{1}{2}-\theta \right) \right] \\*
& = &
\frac{1}{n} \sum_{i=1}^n 2\theta 
\PROB\{ \tilde{h}^{(0)}(Y_1^{i-1},U_i) =0 \} 
  \\*
& \ge &
\frac{1}{n} \sum_{i=1}^n 2\theta 
\PROB\left\{ \sum_{j=1}^{i-1} Y_j < \frac{i-1}{2} \right\} \\*
& = &
\frac{1}{n} \sum_{i=1}^n 2\theta 
\PROB\left\{ \sum_{j=1}^{i-1} (Y_j-\EXP Y_j) < -(i-1)\theta \right\}.
\end{eqnarray*}
Finally, invoke the Berry-Ess\'een  theorem (see, e.g., \cite{ChTe88})
to deduce that there exists a universal constant
$c_2$ such that
$\PROB\left\{ \sum_{j=1}^{i-1} (Y_j-\EXP Y_j) < -(i-1)\theta \right\}
\ge c_2$ for every $2\le i\le n$.)
Thus, even though for every single value
of $\theta$, $\EXP L_1^n(\tilde{h}^{(m)}) - L^*$
converges to zero at a rate of $O(1/n)$, 
the {\em minimax} rate of convergence is, in fact, 
$O(1/\sqrt{n})$.
Since
the upper bound in Theorem \ref{markov} is
independent of the distribution, we see that,
in this sense,
(ignoring logarithmic factors) the order of
magnitude of the bound is the best possible.

\newpage

\section{Prediction with side information}

In this section we apply the same ideas to the seemingly 
more difficult classification (or pattern recognition)
problem. The setup is the following: let
$\{(X_n,Y_n)\}_{-\infty}^{\infty}$
be a stationary and ergodic sequence of pairs
taking values in $\Rd \times \{0,1\}$.
The problem is to predict the value
of $Y_n$ given the data $(X_n,{\cal D}^{n-1})$,
where we denote ${\cal D}^{n-1}=(X_{1}^{n-1},Y_{1}^{n-1})$.
The prediction problem is similar to the one
studied in Section \ref{alg} with the exception
that the sequence of $X_i$'s is also available
to the predictor. One may think about the $X_i$'s
as side information.

We may formalize the prediction problem as follows.
A (randomized) prediction strategy 
is a sequence $g=\{g_i\}_{i=1}^{\infty}$ of decision functions
\[
   g_i: \{0,1\}^{i-1}\times \left(\Rd\right)^i \times [0,1] \to \{0,1\}
\]
so that the prediction formed at time $i$ is
$g_i(y_1^{i-1},x_1^i,U_i)$. 
The {\em normalized cumulative loss} for any fixed
pair of sequences $x_1^n,y_1^n$ is
now
\[
  R_1^n(g,U_1^n) 
   = \frac{1}{n} \sum_{i=1}^n I_{\{g_i(y_1^{i-1},x_1^i,U_i) \neq y_i\}},
\]
We also use the short notation 
%LG kicsereltem itt es mindenhol 
%LG $R_n(g)=R_1^n(g,U_1^n)$.
    $R_1^n(g)=R_1^n(g,U_1^n)$.
Denote the expected loss of the randomized strategy $g$ by
\[
\widehat{R}_1^n(g) = \EXP R_1^n(g,U_1^n).
\]
We assume that the randomizing variables $U_1,U_2,\ldots$
are independent of the process $\{(X_n,Y_n)\}$.

Just like in the case of prediction without side information,
the fundamental limit is given by the Bayes 
probability of error:

\begin{theorem}
\label{bayessideinfo}
For any prediction strategy $g$ and stationary ergodic
process $\{(X_n,Y_n)\}_{n=-\infty}^{\infty}$,
\[
   \liminf_{n\to \infty} R_1^n(g) \ge R^* \quad \mbox{almost surely,}
\]
where 
\[
  R^* = \EXP\left[
\min\left(\PROB\{Y_0=1|Y_{-\infty}^{-1},X_{-\infty}^0\},
   \PROB\{Y_0=0|Y_{-\infty}^{-1},X_{-\infty}^0\}
  \right) \right].
\]
\end{theorem}

The proof of this lower bound is similar to that of Theorem
\ref{bayes}, the details are omitted.
It follows from results of Morvai, Yakowitz, and Gy\"orfi
\cite{MoYaGy96} that there exists a prediction strategy $g$
such that for all ergodic processes, $R_1^n(g)\to R^*$
almost surely. (We omit the details here.)
 The algorithm of Morvai, Yakowitz, and Gy\"orfi,
however, has a very slow rate of convergence
 even for i.i.d.\ processes. 
The main message of this section is a simple 
universal procedure with a practical appeal.
The idea, again, is to combine the decisions
of a small number of simple experts in an
appropriate way.

We define an infinite array
of experts $h^{(k,\ell)}$, $k,\ell=1,2,\ldots$
as follows. 
Let $\P_{\ell}=\{A_{\ell,j}, j=1,2,\ldots,m_{\ell}\}$ be a sequence of 
finite partitions of the feature space $\Rd$, 
and let $G_{\ell}$ be the corresponding quantizer:
\[
G_{\ell}(x)=j, \mbox{ if } x\in A_{\ell,j}.
\]
With some abuse of notation, for any $n$ and
$x_1^n \in \left(\Rd\right)^n$,
we write $G_{\ell}(x_1^n)$ for the sequence 
$G_{\ell}(x_1),\ldots,G_{\ell}(x_n)$.
Fix positive integers $k,\ell$, and for each $s \in \{0,1\}^k$,
$z\in \{1,2,\dots ,m_{\ell}\}^{k+1}$,
and $y \in \{0,1\}$ 
define
\be
{\widehat P}_n^{(k,\ell)}(y,y_1^{n-1},x_1^n,s,z)={
\left|\{k< i < n: y^{i-1}_{i-k}=s, G_{\ell}(x^{i}_{i-k})=z, y_i=y\}
\right|
\over 
\left|\{k < i < n: y^{i-1}_{i-k}=s, 
G_{\ell}(x^{i}_{i-k})=z,\}\right| }, \quad n>k+1.
\ee
$0/0$ is defined to be $1/2$.
Also, for $n \le k+1$ we define 
${\widehat P}_n^{(k,\ell)}(y,y_1^{n-1},x_1^n,s,z)=1/2$.

The expert $h^{(k,\ell)}$ is now defined by
\[
  h^{(k,\ell)}_n(y_1^{n-1},x_1^n) = \left\{ \begin{array}{ll}
   0 & \mbox{if} \ \   
   {\widehat P}_n^{(k,\ell)}
 (0,y_1^{n-1},x_1^n,y_{n-k}^{n-1},G_{\ell}(x_{n-k}^n)) < \frac{1}{2} \\
   1 & \mbox{otherwise,} 
      \end{array}  \right. \quad n=1,2,\ldots
\]
That is, expert $h^{(k,\ell)}$ quantizes the 
sequence $x_1^n$ according to the partition
${\cal P}_{\ell}$, and
looks for all appearances
of the last seen quantized
strings $y_{n-k}^{n-1},G_{\ell}(x_{n-k}^n)$ of length $k$
in the past. Then it predicts according to the larger of the
relative frequencies of 0's and 1's following
the string. 

The proposed algorithm 
combines the predictions of these experts
similarly to that of Section \ref{alg}.
This way both the length of the string to be matched
and the resolution of the quantizer are 
adjusted depending on the data.
The formal definition is as follows:
For any $m=0,1,2,\ldots$, if
 $2^m \le n < 2^{m+1}$, the prediction is based
upon a weighted majority of predictions of the $(2^{m+1})^2$ experts
$h^{(k,\ell)}$, $k,l\le 2^{m+1}$ as follows:
\[
  g_n(y_1^{n-1},x_1^n,u) = \left\{ \begin{array}{ll}
   0 & \mbox{if} \ \ 
\displaystyle{u > 
\frac{\sum_{k,\ell \le 2^{m+1}} 
 h^{(k,\ell)}_n(y_1^{n-1},x_1^n)
 w_n(k,\ell)}
     {\sum_{k,\ell \le 2^{m+1}} w_n(k,\ell)} }  \\
   1 & \mbox{otherwise,}
                     \end{array} \right. 
\] 
where $w_n(k,\ell)$ is the weight of expert $h^{(k,\ell)}$ defined
by the past performance of $h^{(k,\ell)}$ as
\[
w_{2^m}(k,\ell)=1 \quad \mbox{and} \quad  w_n(k,\ell) = 
%LG ugyanaz, mint korabban
  e^{-\eta_m (n-2^m)
 R_{2^m}^{n-1}(h^{(k,\ell)})} \ \mbox{for} \ 2^m < n < 2^{m+1},
\]
where $\eta_m=\sqrt{8\ln(2^{m+1})^2/2^m}$.

\medskip
\noindent
To prove the universality of the method, we need some
natural conditions on
the sequence of partitions. We assume the
following:

\nopagebreak
\medskip
\noindent
(a) the sequence of partitions is nested,
that is, any cell of ${\cal P}_{\ell+1}$ is a
subset of a cell of ${\cal P}_{\ell}$, $\ell=1,2,\ldots$;

\noindent
(b) each partition ${\cal P}_{\ell}$ is finite;

\noindent
(c) if $\diam(A) = \sup_{x,y \in A} \|x-y\|$ denotes the
diameter of a set, then
for each sphere $S$ centered at the
origin
\[
\lim_{\ell\to \infty}
 \max_{j:A_{\ell,j}\cap S\ne \emptyset}\diam(A_{\ell,j})=0.
\]

\medskip
\noindent
{\bf Remark.}
The next theorem states the universality of the
proposed pattern recognition scheme.
The definition of the algorithm is somewhat arbitrary,
we just chose one of the many possibilities.
In this version, at time $n$, only partitions with
indices at most $n$ are taken into account.
It is easy to see that the universality property remains
valid if the number of partitions considered at time $n$ 
is an arbitrary, polynomially increasing function of $n$.
The conditions for the sequence of partitions
again give a lot of liberty to the user. In
applications, the
partitions may be chosen to incorporate some
prior knowledge about the process.
In this paper we merely prove universality
of the scheme. Performance bounds in the
style of Section \ref{mark} for special
types of proceses may be derived, thanks to
the powerful individual sequence bounds.
Here, however, the analysis may be substantially
more complicated.

\begin{theorem}
\label{patt}
Assume that the sequence of 
partitions $\P_{\ell}$ satisfies the three conditions above.
Then the pattern recognition scheme $g$ 
defined above satisfies
\[
 \lim_{n\to \infty} R_1^n(g) = R^* \quad \mbox{almost surely}
\]
for any stationary and ergodic process 
$\{(X_n,Y_n)\}_{n=-\infty}^{\infty}$.
\end{theorem}

\medskip
\noindent
{\bf Proof.}
As in the proof of Theorem \ref{cons}, 
we obtain that for any stationary and ergodic process 
$\{(X_n,Y_n)\}_{n=-\infty}^{\infty}$,
\begin{eqnarray*}
\limsup_{n\to \infty} R_1^n(g,U_1^n)
& \le & \limsup_{n\to \infty}
\left( 
\min_{\scriptsize
{\begin{array}{c} k=1,2,\ldots \\ \ell=1,2,\ldots,n-1 \end{array}}}
 R_1^n(h^{(k,\ell)}) + 
2c \sqrt{\frac{\log_2 n + 1}{n}} \right)      \\*
& = &
\limsup_{n\to \infty} 
\min_{\scriptsize
{\begin{array}{c} k=1,2,\ldots \\ \ell=1,2,\ldots,n-1 \end{array}}}
 R_1^n(h^{(k,\ell)})
 \quad \mbox{almost surely}.
\end{eqnarray*}
Thus, it remains to show that
\[
\limsup_{n\to \infty}
\min_{\scriptsize
{\begin{array}{c} k=1,2,\ldots \\ \ell=1,2,\ldots,n-1 \end{array}}}
 R_1^n(h^{(k,\ell)}) \le R^*  \quad \mbox{almost surely}.
\]
To prove this, we use the
following lemma, whose proof is easily obtained by
copying that of Lemma \ref{fststat}:

\begin{lemma}
\label{stat}
For each $k,\ell \ge 1$, there exists a positive number
$\epsilon_{k,\ell}$ such that for any fixed $\ell$,
$\lim_{k\to \infty} \epsilon_{k,\ell} =0$ and 
\[
\limsup_{n\to\infty} R_1^n(h^{(k,\ell)})
  \le R^*_{(\ell)} + \epsilon_{k,\ell},
\]
where 
\[
R^*_{(\ell)} = \EXP\left[ 
\min\left( \PROB\{Y_0=1|Y_{-\infty}^{-1},G_{\ell}(X^{0}_{-\infty})\},
\PROB\{Y_0=0|Y_{-\infty}^{-1},G_{\ell}(X^{0}_{-\infty})\}\right)
  \right].
\]
\end{lemma}

\medskip
\noindent
Now we return to the proof of Theorem \ref{patt}.
Since the sequence of partitions $\P_{\ell}$ is nested, and by (c),  
the sequences 
\[\PROB\{Y_0=1|Y_{-\infty}^{-1},G_{\ell}(X^{0}_{-\infty})\}
\quad \mbox{and} 
\quad\PROB\{Y_0=0|Y_{-\infty}^{-1},G_{\ell}(X^{0}_{-\infty})\}
\quad l=1,2,\ldots
\]
are martingales and they converge almost surely to 
\[\PROB\{Y_0=1|Y_{-\infty}^{-1},X^{0}_{-\infty}\}
\quad \mbox{and} 
\quad\PROB\{Y_0=0|Y_{-\infty}^{-1},X^{0}_{-\infty}\}.
\]
Thus, it follows from Lebesgue's dominated convergence theorem
that
\[
\lim_{l\to \infty} R^*_{(\ell)}
 = \EXP\left[ 
\min\left( \PROB\{Y_0=1|Y_{-\infty}^{-1},X^{0}_{-\infty}\},
\PROB\{Y_0=0|Y_{-\infty}^{-1},X^{0}_{-\infty}\}\right)
  \right]=R^*.
\]
Now it follows easily that
\[
\limsup_{n\to \infty}
\min_{\scriptsize
{\begin{array}{c} k=1,2,\ldots \\ \ell=1,2,\ldots,n-1 \end{array}}}
 R_1^n(h^{(k,\ell)}) \le R^*  \quad \mbox{almost surely},
\]
and the proof of the theorem is finished. \qed

\newpage

\section{Appendix}

Here we describe two results which are used in the analysis.
The first is due to Breiman \cite{Bre60}, and its proof may
also be found in Algoet \cite{Alg94}.

\begin{lemma} \label{Breimanergod}
{\sc Breiman's generalized ergodic theorem \cite{Bre60}.}
Let $Z=\{Z_i\}^{\infty}_{-\infty}$  be a stationary and ergodic time
series.   Let $T$ denote the left shift operator. Let $f_i$ be a
sequence of real-valued functions such that for some function $f$,
$f_i(Z)\to f(Z)$ almost surely.  
Assume that $\EXP\sup_i|f_i(Z)|<\infty$. Then
$$
\lim_{t\to\infty} {1\over n} \sum_{i=1}^{n} f_i(T^i Z)=\EXP f(Z)
$$
almost surely. 
\end{lemma}

The second is the Hoeffding-Azuma inequality
for sums of bounded martingale differences:

\begin{lemma}
\label{azuma}
{\sc Hoeffding \cite{Hoe63}, Azuma \cite{Azu67}.} \ \
Let $X_1,X_2,\ldots$ be a sequence of random variables,
and assume that $V_1,V_2,\ldots$ is a martingale difference
sequence with respect to $X_1,X_2,\ldots$. Assume furthermore that
there exist random variables $Z_1,Z_2,\ldots$
and nonnegative constants $c_1,c_2,\ldots$  such that
for every $i>0$ $Z_i$ is a function of $X_1,\ldots,X_{i-1}$,
and
\[
    Z_i \le V_i \le Z_i+ c_i \; \; \mbox{with probability one.}
\]
Then for any $\epsilon >0$ and $n$
\[
   \PROB\left\{ \sum_{i=1}^n V_i \ge \epsilon \right\}\le
      e^{-2\epsilon^2/ \sum_{i=1}^n c_i^2}
\]
and
\[
   \PROB\left\{ \sum_{i=1}^n V_i \le -\epsilon \right\}\le
      e^{-2\epsilon^2/ \sum_{i=1}^n c_i^2}.
\]
\end{lemma}

\medskip
\noindent
{\bf Acknowledgement.} We thank Nicol\'o Cesa-Bianchi  
for teaching us all wee needed to know about prediction with
expert advise. We are also grateful to Sid Yakowitz for
illuminating discussions and to the referees for 
a very careful reading of the manuscript and for
valuable suggestions.
We also thank M\'arta Horv\'ath for useful conversations.

%\bibliography{book1,book2,sajat,cikkek}

\end{document}